\documentclass{article}

\usepackage{amssymb}

\newcommand{\real}{{\mathbb R}}
\newcommand{\iin}{\!\in\!}
\newcommand{\Lip}{\mbox{\rm Lip}}
\newcommand{\Ub}{{\bf U_b}}
\newcommand{\Meas}{{\bf M}}
\newcommand{\UMeas}{{\bf M_u}}
\newcommand{\CMeas}{{\bf M_\sigma}}
\newcommand{\CUMeas}{{\bf M_{u \sigma}}}
\newcommand{\compln}[1]{\widehat{#1}}
\newcommand{\Coz}{\mbox{\bf Coz}}
\newcommand{\Open}[1]{{\cal O}({#1})}
\newcommand{\wstar}{weak\raisebox{1mm}{$\ast$}\ }
\newcommand{\bigsep}{\mbox{\Large $|$}}
\newcommand{\om}{\widetilde{m}}
\newcommand{\qed}{\hfill $\Box$\vspace{3mm}}

\newtheorem{theorem}{Theorem}
\newtheorem{lemma}[theorem]{Lemma}
\newtheorem{corollary}[theorem]{Corollary}

\title{Uniform measures and countably additive measures}

\author{Jan Pachl  \\ Toronto, Ontario, Canada}
\date{April 6, 2007}

\begin{document}
\maketitle

\begin{abstract}
Uniform measures are defined as the functionals on the space of bounded
uniformly continuous functions that are continuous on
bounded uniformly equicontinuous sets.
If every cardinal has measure zero then every countably
additive measure is a uniform measure.
The functionals sequentially continuous on bounded uniformly equicontinuous sets
are exactly uniform measures on the separable modification
of the underlying uniform space.
\end{abstract}


\section{Introduction}
    \label{section:introduction}

The functionals that we now call uniform measures
were originally studied by
Berezanski\v{\i}~\cite{Berezanskii1968},
Csisz\'{a}r~\cite{Csiszar1971},
Fedorova~\cite{Fedorova1967}
and LeCam~\cite{LeCam1970}.
The theory was later developed in several directions by a number of other authors;
see the references in~\cite{Pachl2006} and~\cite{Rice1981}.

Uniform measures need not be countably additive,
but they have a number of properties that have traditionally been
formulated and proved for countably additive measures,
or countably additive functionals on function spaces.
The main result in this paper,
in section~\ref{section:dspaces}, is that countably additive measures
are uniform measures on a large class of uniform spaces
(on all uniform spaces, if every cardinal has measure zero).

Section~\ref{section:cuniformmeasures} deals with the functionals
that behave like uniform measures on sequences of functions;
or, equivalently, like countably additive measures on bounded
uniformly equicontinuous sets.
In the case of a topological group with its right uniformity,
these functionals were defined by Ferri and Neufang~\cite{Ferri-Neufang2006}
and used in their study of topological centres in convolution algebras.

\section{Notation}
    \label{section:notation}

In the whole paper, linear spaces are assumed to be over the field $\real$ of reals.
Uniform spaces are assumed to be Hausdorff.
Uniform spaces are described by uniformly continuous pseudometrics (\cite{Gillman-Jerison1960}, Chap.~15),
abbreviated u.c.p.

When $d$ is a pseudometric on a set $X$, define
\[
\Lip ( d )  =  \{ f: X \rightarrow \real \; \bigsep \; |f(x)| \leq 1 \; \mbox{\rm and} \;
| f(x) - f(x') | \leq d(x,x') \; \mbox{\rm for all} \; x,x' \in X  \} \; .
\]
Then $\Lip(d)$ is compact in the topology of pointwise convergence on $X$,
as a topological subspace of the product space $ \real^X $.

When $X$ is a uniform space, denote by $\Ub(X)$ the space of bounded uniformly continuous
functions $ f : X \rightarrow \real $ with the norm
$
\| f \| = \sup \{ \; | f(x) | \;
    | \; x \in X \}
$.
Let $ \Coz(X) $ be the set of all cozero sets in $X$;
that is, sets of the form
$ \{ x\in X \; | \; f(x) \neq 0 \} $
where $ f \in \Ub (X) $.
Let $ \sigma ( \Coz(X) ) $ be the sigma-algebra of subsets of $X$
generated by $ \Coz(X) $.

When $d$ is a pseudometric on a set $X$, denote by
$ \Open{d} $ the collection of open sets in the (not necessarily Hausdorff)
topology defined by $d$.
Note that if $d$ is a u.c.p. on a uniform space $X$ then
$ \Open{d} \subseteq \Coz (X) $.

Denote by $\Meas(X)$ the norm dual of $\Ub(X)$, and consider three subspaces of $ \Meas(X) $:
\begin{enumerate}
\item
$ \UMeas(X) $ is the space of those $ \mu \in \Meas(X) $ that are continuous on
$\Lip(d)$ for every u.c.p. $d$ on $X$,
where $\Lip(d)$ is considered with the topology of pointwise convergence on $X$.
The elements of $ \UMeas(X) $ are called {\em uniform measures on $X$\/}.
\item
$ \CMeas(X) $ is the space of $ \mu \in \Meas(X) $ for which
there is a bounded (signed) countably additive measure
$m$ on the sigma-algebra $ \sigma ( \Coz(X) ) $ such that
\[
\mu (f) = \int f \mbox{\rm d} m \;\; \mbox{\rm for} \;\; f \in \Ub(X) \; .
\]
\item
$ \CUMeas(X) $ is the space of those $ \mu \in \Meas(X) $
that are sequentially continuous on $ \Lip(d) $ for each u.c.p. $d$.
That is, $ \lim_n \mu ( f_n ) = 0 $ whenever
$d$ is a u.c.p. on $X$, $ f_n \in \Lip(d) $ for $ n = 1, 2, \ldots $,
and $ \lim_n f_n (x) = 0 $ for each $ x \in X $.
\end{enumerate}

When $X$ is a topological group $G$ with its right uniformity,
$ \CUMeas(X) $ is the space $ \mbox{\rm Leb}^s (G) $ in the notation of~\cite{Ferri-Neufang2006}.

Clearly $ \UMeas(X) \subseteq \CUMeas(X) $ for every uniform space $X$.
By Lebesgue's dominated convergence theorem (\cite{FremlinMT1}, 123C),
$ \CMeas(X) \subseteq \CUMeas(X) $ for every $X$.

For any uniform space $X$, let $ cX $ be the set $X$ with the weak uniformity induced
by all uniformly continuous functions from $X$ to $\real$ (\cite{Isbell1964}, p.~129).
Let $eX$ be the {\em cardinal reflection\/} $X_{\aleph_1}$ (\cite{Isbell1964}, p.~52 and~129),
also known as the {\em separable modification\/} of $X$.
Thus $ eX $ is a uniform space on the same set as $X$, and a pseudometric on $X$
is a u.c.p. on $eX$ if and only if it is a separable u.c.p. on $X$.
Note that $ \Ub (X) = \Ub (eX) = \Ub(cX) $ and
$ \Meas(X) = \Meas(eX) = \Meas(cX) $.

Let $\aleph$ be a cardinal number, and let $A$ be a set of cardinality $\aleph$.
As in~\cite{Granirer1967},
say that $ \aleph $ has {\em measure zero\/} if
$ m(A) = 0 $ for every non-negative countably additive measure
$m$ defined on the sigma-algebra of all subsets of $A$ and
such that $ m( \{ a \} ) = 0 $ for all $ a \in A $.
A related notion, not used in this paper,
is that of a {\em nonmeasurable cardinal\/}
as defined by Isbell~\cite{Isbell1964},
using two-valued measures $m$ in the preceding definition.

It is not known whether every cardinal has measure zero.
The statement that every cardinal has measure zero is consistent with the usual
axioms of set theory.
A detailed discussion of this and related properties of cardinal numbers
can be found in~\cite{FremlinMT5} and~\cite{Jech1997}.

Let $d$ be a pseudometric on a set $X$.
A collection $ \cal W $ of nonempty subsets of $X$ is {\em uniformly $d$-discrete\/}
if there exists $ \varepsilon > 0$ such that $ d(x,x') \geq \varepsilon $ whenever
$ x \!\in\! V, x' \!\in\! V' $, $ V, V' \!\in\! \cal W $, $ V \neq V' $.
A set $ Y \subseteq X $ is {\em uniformly $d$-discrete\/} if the collection of
singletons $ \{ \{ y \} \; | \; y \in Y \} $ is uniformly $d$-discrete.

Let $X$ be a uniform space.
A set $ Y \subseteq X $ is {\em uniformly discrete\/} if there exists
a u.c.p. $d$ on $X$ such that $Y$ is uniformly $d$-discrete.
Say that $X$ is a {\em (uniform) D-space\/}~\cite{Pachl1976}
if the cardinality of every uniformly discrete subset of $X$
has measure zero.

This generalizes the notion of a topological D-space as defined by Granirer~\cite{Granirer1967}
and further discussed by Kirk~\cite{Kirk1973} in the context of
topological measure theory.
A topological space $T$ is a D-space in the sense of~\cite{Granirer1967}
if and only if $T$ with its fine uniformity (\cite{Isbell1964}, I.20) is a uniform D-space.
If $X$ is a uniform space and $ Y \subseteq X $ is uniformly discrete in $X$
then $Y$ is also uniformly discrete in $X$ with its fine uniformity.
Therefore, if $X$ is a topological D-space in the sense of~\cite{Granirer1967}
then it is also a uniform D-space.

Since the countable infinite cardinal $ \aleph_0 $ has measure zero,
every uniform space $X$ such that $ X = eX $ is a D-space.
Thus every uniform subspace of a product of separable metric spaces is a D-space.
Moreover, the statement that {\em every\/} uniform space is a D-space is consistent with the usual
axioms of set theory.


\section{Measures on uniform D-spaces}
    \label{section:dspaces}

The uniform spaces $X$ for which $ \UMeas(X) \subseteq \CMeas(X) $
were investigated by several authors
\cite{Berezanskii1968}
\cite{Fedorova1967} \cite{Fedorova1974}
\cite{Frolik1975} \cite{LeCam1970}.
The opposite inclusion $ \CMeas(X) \subseteq \UMeas(X) $ has not attracted as much attention.
Theorem~\ref{th:Dspace} in this section characterizes the uniform spaces $X$ for which
$ \CMeas(X) \subseteq \UMeas(X) $.

\begin{lemma}
     \label{lemma:sigmadiscrete}
Let $d$ be a pseudometric on a set $X$, and $ \varepsilon > 0 $.
Then there exist sets $ {\cal W}_n $ of nonempty subsets of $X$, $ n = 1,2, \ldots $,
such that
\begin{enumerate}
\item
$ \bigcup_{n=1}^\infty {\cal W}_n $ is a cover of $X$;
\item
for each $n$, $ {\cal W}_n \subseteq \Open{d} $;
\item
for each $n$, the $d$-diameter of each $ V \in {\cal W}_n $ is at most $ \varepsilon $;
\item
each $ {\cal W}_n $ is uniformly $d$-discrete.
\end{enumerate}
\end{lemma}

\noindent
The lemma is essentially the theorem of A.H.~Stone about $\sigma$-discrete covers
in metric spaces.
For the proof, see the proof of~4.21 in~\cite{Kelley1967}.

The next theorem is the main result of this paper.
It generalizes a known result about separable measures on completely regular
topological spaces --- Proposition 3.4 in~\cite{Kirk1973}.

\begin{theorem}
    \label{th:Dspace}
    For any uniform space $X$, the following statements are equivalent:
    \begin{description}
    \item[{\it(i)}] $X$ is a uniform D-space.
    \item[{\it(ii)}] $ \CMeas(X) \subseteq \UMeas(X) $.
    \end{description}
\end{theorem}

In view of Theorem~\ref{th:Dspace} and the remarks in section~\ref{section:notation},
the statement that $ \CMeas(X) \subseteq \UMeas(X) $ for every uniform space $X$
is consistent with the usual axioms of set theory.

\noindent
{\bf Proof.}
This proof is adapted from the author's
unpublished manuscript~\cite{Pachl1976}.

To prove that (i) implies (ii), let $X$ be a D-space.
To show that $ \CMeas(X) \subseteq \UMeas(X) $,
it is enough to show that $ \mu \in \UMeas(X) $ for every non-negative
$ \mu \in \CMeas(X) $,
in view of the Jordan decomposition of countably additive measures
(\cite{FremlinMT2}, 231F).
Take any $ \mu \in \CMeas(X) $, $ \mu \geq 0$ and any $ \varepsilon > 0 $.
Let $m$ be the non-negative countably additive measure on $ \sigma( \Coz(X) ) $
such that
$
\mu (f) = \int f \mbox{\rm d} m
$
for $ f \in \Ub(X) $.

Let $d$ be a u.c.p. on $X$, and  $ \{ f_\alpha \}_\alpha $
a net of functions $ f_\alpha \in \Lip(d) $
such that $ \lim_\alpha f_\alpha ( x ) = 0 $ for every $ x \in X $.
Our goal is to prove that $ \lim_\alpha \mu( f_\alpha ) = 0 $.

For the given $X$, $d$ and $\varepsilon$,
let $ {\cal W}_n $ be as in Lemma~\ref{lemma:sigmadiscrete}.
If $ V \in {\cal W}_n $ for some $n$ then choose a point $ x_V \in V $.
Let $ T_n = \{ x_V \; | \; V \in {\cal W}_n \} $ for $n=1,2, \ldots$.

Fix $n$ for a moment.
For each subset $ {\cal W'} \subseteq {\cal W}_n $
we have $ \bigcup {\cal W'} \in \Open{d} \subseteq \Coz (X) $.
Thus for each $ S \subseteq T_n $ we may define
$ \om ( S ) = m ( \, \bigcup \{ V \in {\cal W}_n \; | \; x_V \in S \} \, )$,
and $\om$ is a countably additive measure defined on all subsets of $T_n$.
Since the set $ T_n $ is uniformly discrete and $X$ is a D-space,
it follows that the cardinality of $ T_n $ is of measure zero,
and there exists a countable set $ S_n \subseteq T_n $ such that
\[
m ( \, \bigcup \{ V \in {\cal W}_n \; | \; x_V \in T_n \setminus S_n \} \, )
= \om ( T_n \setminus S_n ) = 0 .
\]

Denote $ P = \bigcup_{n=1}^\infty S_n $ and
$ Y = \{ x \in X \; | \; d( x , P ) \leq \varepsilon \} $.
If $ V \in {\cal W}_n $ for some $n$ and $ x_V \in P $ then $ V \subseteq Y $,
by property 3 in Lemma~\ref{lemma:sigmadiscrete}.
Therefore
\[
X \setminus Y \subseteq \bigcup_{n=1}^\infty \bigcup \{ V \in {\cal W}_n \; | \; x_V \in T_n \setminus S_n \}
\]
and $ m ( X \setminus Y ) = 0 $.

Define
$ g_\alpha (x) = \sup_{\beta \geq \alpha} | f_\beta (x) | $
for $ x \in X $.
Then $ g_\alpha \in \Lip(d) $,
$ g_\alpha \geq g_\beta $ for $ \alpha \leq \beta $, and
$ \lim_\alpha g_\alpha ( x ) = 0 $ for every $ x \in X $.

Since the set $P$ is countable, there is an increasing sequence
of indices $ \alpha (n) $, $ n=1,2, \ldots$, such that
$ \lim_n g_{\alpha(n)} ( x ) = 0 $ for every $ x \in P $,
hence $ \lim_n g_{\alpha(n)} ( x ) \leq \varepsilon $ for every $ x \in Y $.
Thus
\[
\lim_\alpha | \mu( f_\alpha ) |
\leq \lim_\alpha \mu( g_\alpha )
\leq \lim_n \mu( g_{\alpha(n)} )
= \lim_n \left( \int_Y g_{\alpha(n)} {\rm d}m
+ \int_{X \setminus Y} g_{\alpha(n)} {\rm d}m \right)
\leq \varepsilon \, m(X)
\]
which proves that $ \lim_\alpha \mu( f_\alpha ) = 0 $.

To prove that (ii) implies (i), assume that $X$ is not a D-space.
Thus there is a u.c.p. $d$ on $X$, a subset $ P \subseteq X $
and a non-negative countably additive measure $m$ defined on all subsets of $P$ such that
\begin{itemize}
\item
$ d(x,y) \geq 1 $ for $ x,y \in P $, $ x\neq y $;
\item
$ m (x) = 0 $ for each $ x \in P $;
\item
$ m ( P ) = 1 $.
\end{itemize}
Define $ \mu (f) = \int_P f {\rm d}m $ for $ f \in \Ub(X)$.
Clearly $ \mu \in \CMeas(X) $.

For any set $ S \subseteq P $, define
the function $ f_S \in \Lip(d) $ by
$ f_S (x) = \min ( 1 , d(x,S) ) $ for $ x \in X$.
Then $ f_S (x) = 0 $ for $ x \iin S $ and
$ f_S (x) = 1 $ for $ x \iin P \setminus S $.
Let $\cal F$ be the directed set of all finite subsets of $P$ ordered by inclusion.
We have $ \lim_{S\in {\cal F}} f_S (x) = f_P (x) $ for each $ x \in X $,
$ \mu ( f_S ) = 1 $ for every $ S \in {\cal F} $,
and $ \mu ( f_P ) = 0 $.
Thus $ \mu \not\in \UMeas(X) $.
\qed

The inclusion $ \CMeas(X) \subseteq \UMeas(cX) $ in the following corollary
is Theorem~2.1 in~\cite{Fedorova1982}.

\begin{corollary}
    \label{cor:cmeasex}
If $X$ is any uniform space then
$ \CMeas(X) \subseteq \UMeas(eX) \subseteq \UMeas(cX) $.
\end{corollary}

\noindent
{\bf Proof.}
As is noted above, $ eX $ is a D-space for any $X$.
Thus $ \CMeas(eX) \subseteq \UMeas(eX) $ by Theorem~\ref{th:Dspace}.
From the definitions of $\CMeas(X)$, $eX$ and $cX$ we get $ \CMeas(X) = \CMeas(eX) $ and
$ \UMeas(eX) \subseteq \UMeas(cX) $.
\qed

Corollary~\ref{cor:cmeasex} follows also from Theorem~\ref{th:cumeasures} in the next section:
$ \CMeas(X) \subseteq \CUMeas(X) = \UMeas(eX) $.


\section{Countably uniform measures}
    \label{section:cuniformmeasures}

In this section we compare the spaces $ \CUMeas(X) $ and $ \UMeas(X) $.

\begin{theorem}
    \label{th:cumeasures}
If $X$ is any uniform space then $ \CUMeas(X) = \UMeas(eX) $.
\end{theorem}

\noindent
{\bf Proof.}
To prove that $ \CUMeas(X) \subseteq \UMeas(eX) $,
note that if a pseudometric $d$ is separable then
$ \Lip(d) $ with the topology of pointwise convergence is metrizable,
and therefore sequential continuity on $ \Lip(d) $ implies continuity.

To prove that $ \UMeas(eX) \subseteq \CUMeas(X) $,
take any $ \mu \in \UMeas(eX) $.
Let $d$ be a u.c.p. on $X$, $ f_n \in \Lip(d) $ for $ n = 1, 2, \ldots $,
and $ \lim_n f_n (x) = 0 $ for each $ x \in X $.
Define a pseudometric $\widetilde{d}$ on $X$ by
\[
\widetilde{d}(x,y) = \sup_n | f_n (x) - f_n (y) | \;\;\; \mbox{\rm for} \;\;\; x,y \in X \; .
\]
Then $\widetilde{d}$ is a separable u.c.p. on $X$, hence a u.c.p. on $eX$,
and $ f_n \in \Lip(\widetilde{d}\, ) $ for $ n = 1, 2, \ldots $.
Therefore $ \lim_n \mu ( f_n ) = 0 $.
\qed

In view of Theorem~\ref{th:cumeasures}, spaces $ \CUMeas(X) $ have all the properties
of general $ \UMeas(X) $ spaces.
For example, every $ \UMeas(X) $ is \wstar sequentially complete~\cite{Pachl2006},
and the positive part $\mu^+$ of every $ \mu \in \UMeas(X) $ is in $\UMeas(X)$
\cite{Berezanskii1968} \cite{Fedorova1967} \cite{LeCam1970}.
Therefore the same is true for $ \CUMeas(X) $.

By Theorem~\ref{th:cumeasures}, if $ X = eX $ then $ \UMeas(X) = \CUMeas(X) $
(cf.~\cite{Ferri-Neufang2006}, 2.5(iii)).
To see that the equality $ \UMeas(X) = \CUMeas(X) $ does not hold in general,
first consider a uniform space $X$ that is not a uniform D-space.
Since $ \CMeas(X) \subseteq \CUMeas(X) $, from Theorem~\ref{th:Dspace} we get
$ \UMeas(X) \neq \CUMeas(X) $.
However, that furnishes an actual counterexample only if there exists a cardinal
that is not of measure zero.
Next we shall see that, even without assuming the existence of such a cardinal,
there is a space $X$ such that $ \UMeas(X) \neq \CUMeas(X) $.

Let $ \compln{X} $ denote the completion of a uniform space $X$.
Pelant~\cite{Pelant1975} constructed a complete uniform space X
for which $ eX $ is not complete.
For such $X$,
there exists an element $ x \in \compln{eX} \setminus X $.
Every $ f \iin \Ub(X) = \Ub(eX) $ uniquely extends to $ \compln{f} \iin \Ub ( \compln{eX} )$.
Let $ \delta_x \iin \Meas(X) $ be the Dirac measure at $x$;
that is, $ \delta_x (f) = \compln{f} (x) $ for $ f \iin \Ub(X) $.
Then $ \delta_x \iin \UMeas(eX) $, therefore
$ \delta_x \iin \CUMeas(X) $ by Theorem~\ref{th:cumeasures}.
On the other hand, $ \delta_x \not\in \UMeas(X) $,
since $ \delta_x $ is a multiplicative functional on $ \Ub(X) $
and $ x \not\in \compln{X} $ (\cite{Pachl2006}, section~6).
Thus $ \UMeas(X) \neq \CUMeas(X) $.


\end{document}